\input amstex
\documentstyle{amsppt}
\magnification=\magstep1
\parskip=6pt
\baselineskip=14pt
\pagewidth{5.0in}
\pageheight{7.0in}
\CenteredTagsOnSplits
\TagsOnRight
\define\real{{\Bbb R}}
\magnification=\magstep1
\parskip=6pt
\def \gl3p {{$GL(3)^+$}}
\define\comp{{\Bbb C}}
\def\cala{{\Cal A}}
\def\calb{{\Cal B}}
\def\calg{{\Cal G}}
\def\calm{{\Cal M}}

\def\qq{{\Bbb Q}}
\define\zz{{\Bbb Z}}

\def\inc{\smash{\mathop{\hookrightarrow}\limits^i}}
\def\mapright#1{\smash{\mathop{\longrightarrow}\limits^{#1}}}
\def\mapdown#1{\Big\downarrow\rlap{$ \vcenter{\hbox{$ \scriptstyle#1$}}$}}

\topmatter
\title A Simple Geometric Representative for $\mu$ of a Point
\endtitle
\author  Lorenzo Sadun  \endauthor
\thanks{Department of Mathematics,
University of Texas, Austin, Texas 78712.
\hfill \break \indent This work is partially supported by
an NSF Mathematical Sciences Postdoctoral Fellowship and
Texas Advanced Research Project grant ARP-037. 
\hfill \break \indent Email: sadun\@math.utexas.edu, \qquad FAX: 512-471-9038
\hfill \break \indent Running header: Simple Geometric Representative for
$\mu$ of a Point
}
\endthanks
\abstract
{For $SU(2)$ (or $SO(3)$) Donaldson theory on a 4-manifold $X$,
we construct a simple geometric representative for $\mu$ of a point. 
Let $p$ be a generic point in $X$.
Then the set $\{ [A] | F_A^-(p) $ is reducible $\}$, with coefficient
-1/4 and appropriate orientation, is our desired geometric representative.
The construction is an exercise in real algebraic geometry in the style
of Ehresmann and Pontryagin.}
\endabstract
\endtopmatter
\document

\noindent {\bf 1. Background and Statement of Results}

In the past decade, an industry has developed 
studying the homology of moduli spaces, thereby shedding light on 
the topology or geometry of underlying manifolds.  
The best known example is Donaldson's work on gauge
theory in 4 dimensions [DK]. Donaldson's polynomial invariants measure the
fundamental classes of moduli spaces of 
anti-self-dual connections over an orientable 4-manifold, 
giving information about the 
differentiable structure of that manifold. 

Let $X$ be an oriented 4-manifold, let $G=SU(2)$ or $SO(3)$
and let $\calb_k$ be the space of 
connections (up to gauge equivalence) 
on $P_k$, the principal $G$ bundle of instanton
number $k$ over $X$.  Let $\calb^*_k$ (resp. $\tilde \calb^*_k$) be the 
space of irreducible connections, (resp. irreducible framed connections) on $P_k$, modulo gauge equivalence.
$\tilde \calb^*_k$ is a principal $SO(3)$ bundle over $\calb^*_k$.

Donaldson [D1, D2] defined a map $\tilde \mu: H_i(X,\qq) \to
H^{4-i}(\tilde \calb^*_k,\qq)$, $i=$1, 2, 3, whose image freely generates 
the rational cohomology of $\tilde \calb^*_k$.  
For $\Sigma$ a 1, 2, or 3-cycle in $X$, the class $\tilde \mu([\Sigma])$
descends to a cohomology class on $\calb^*_k$, which is then denoted
$\mu([\Sigma])$.  The classes $\mu([\Sigma])$, together with an additional 
4-dimensional class, freely generate the cohomology of $\calb^*_k$.
The additional class can be viewed as $\mu$ of the point class
$[x] \in H_0(X)$.  In this view, $\mu$ maps $H_i(X)$ to $H^{4-i}(\calb^*_k)$,
where $i$ now ranges from 0 to 3, and the image of the $\mu$ map
freely generates $H^*(\calb^*_k,\qq)$.

This gives a polynomial invariant on the homology of $X$, 
the action of $\mu$ of the
elements of $H_*$ on the ``fundamental class'' of $\calm_k$.  Formally, for
elements $[\Sigma_1],\ldots,[\Sigma_n] \in H_*(X)$,
we write
$$q([\Sigma_1],\ldots,[\Sigma_n]) = \mu([\Sigma_1]) \smile \cdots \smile \mu([\Sigma_n])
[\calm_k]. \eqno(1) $$
The ``fundamental class of $\calm_k$'' is usually not well defined, as 
$\calm_k$ is typically not compact.  To make sense of (1) one must
compactify $\calm_k$ and show that the classes $\mu([\Sigma])$ extend properly
to the compactification of $\calm_k$.  This is usually done with geometric
representatives.
One finds finite-codimension varieties $V_\Sigma$ in $\calb$ that are,
roughly speaking, Poincare dual to $\mu([\Sigma])$.  One then attempts to
count points in $V_{\Sigma_1} \cap \cdots \cap V_{\Sigma_n} \cap \calm_k$.
To make a topological invariant one must show that the intersection
points can be bounded away from the ends of $\calm_k$. This requires 
careful analysis of the
bubbling-off phenomena that make $\calm_k$ noncompact.

The success of such a program can depend on good choices of geometric
representatives.  For example, for 2-dimensional Yang-Mills theory, the
generalized Newstead conjecture resisted abstract analysis 
until Weitsman [We] found
a set of simple geometric representatives for the problem.  Using these,
it was fairly easy to characterize the points in $\cap_i V_{\Sigma_i} \cap \calm$, compute the invariants, and prove the conjecture.

For Donaldson theory, fairly simple geometric representatives 
have been found for the 1, 2, and 3-dimensional classes.
In each case, the geometric representative of
$\mu([\Sigma])$ is the set of connections that satisfy a 
simple condition when restricted to $\Sigma$.  
Until now, however, there has not been any similar description
of $\mu([p])$, where $p$ is a single point, 
in terms of data at that point.
The purpose of this paper is to provide such a description.  
For any point $p \in X$, 
let $\nu_p = \{ [A] \in \calb^*_k | F_A^-$ is reducible at $p \}$.  Here 
$F_A^- = (F_A - *F_A)/2$ is the anti-self-dual part of the curvature $F_A$, 
and by ``reducible at $p$'' we mean that the components $F_{ij}^-(p)$
are all colinear as elements of the Lie algebra of $G$.  The main
theorem is

\noindent {\bf Theorem 1:} {\it $\nu_p$ is a geometric representative  
of $-4 \mu([p])$.}

The proof proceeds in stages.  In section 2, we review some classical
real algebraic geometry
and construct a simple representative of the first Pontryagin class
$p_1$ of canonical 
$SO(3)$ bundles over Grassmannians of real oriented 3-planes.  
The construction is essentially due to Pontryagin [P] and 
Ehresmann [E], but their techniques seem to have been generally forgotten.
In section 3, we extend this analysis to $BSO(3)$
and construct an explicit isomorphism between
a space of connections on a neighborhood of 
the point $p$ and $ESO(3)$.  Pulling the representative of $p_1(ESO(3))$
back by this isomorphism gives $\nu_p$, and fixes the orientation. 

To be useful for Donaldson theory, $\nu_p$ must be transverse to
the moduli spaces $\calm_k$ and extend to the compactification
of $\calm_k$.  These issues are discussed in section 4, where we also 
discuss a possible topological application of this representative.

\medskip

\noindent {\bf 2. Cohomology of Real Grassmannians}

\medskip

Let $V_N$ be the space of real, rank 3, $3 \times N$ matrices.  
Equivalently, $V_N$ is the Stiefel manifold of triples of linearly
independent vectors in $\real^N$.  Let $V^0_N$ be the triples
of orthonormal vectors in $\real^N$.  The group $SO(3)$ acts 
freely on both
spaces by left multiplication.  Let $B_N$ be the quotient of $V_N$ by 
$SO(3)$ and let $G_N$ be the quotient of $V^0_N$ by $SO(3)$.
$G_N$ is the Grassmannian of oriented 3-planes in $\real^N$.
We will denote by $\pi$ both natural projections, from $V_N$ to
$B_N$ and from $V^0_N$ to $G_N$.  The Gram-Schmidt process gives
a natural bundle map from $V_N$ to $V^0_N$, which we denote $\rho$.
$\rho$  itself defines trivial $\real^6$ bundles  
$V_N \to V^0_N$ and $B_N \to G_N$.
Inclusion of  $V^0_N$ in $V_N$ defines a natural section. 
In short, we have the commutative diagram
$$ \matrix 
V_N & \mapright{\rho} & V^0_N \cr
\mapdown{\pi} & & \mapdown{\pi} \cr
B_N & \mapright{\rho} & G_N 
\endmatrix
\eqno(2) $$
$B_N$ and $G_N$ have the same topology. 
 
\noindent {\bf Theorem 2.} {\it 
Let $\nu_N=\{ m \in V_N | $ first 3 columns of $m$ have rank $\le 1 \}$.  
Then $\pi(\nu_N)$ is Poincare dual to a generator of $H^4(B_N)$.
By choosing orientations correctly, this generator may be taken to
be the first Pontryagin class of the bundle $V_N \to B_N$.} 

\noindent {\it Proof:} The proof is an application 
of some general computations of Pontryagin [P] and Ehresmann [E].
(Indeed, theorem 2 was almost certainly known to Pontryagin).  
Within the 9 dimensional space of real $3 \times 3$ matrices, 
the rank $\le 1$ matrices form a closed codimension-4 set.
$\pi(\nu_N)$ is thus a closed codimension-4
submanifold of $B_N$, and so is dual to some (possibly zero)
element of $H^4$.  We construct a generator of $H_4(B_N)$ and show 
it intersects $\pi(\nu_N)$ exactly once, establishing that 
$\pi(\nu_N)$ is a generator of $H^4$.  The sign, relative to $p_1$, is 
determined separately.

We begin with a cell decomposition of $G_N$.  Consider the set of
$3 \times N$ matrices of the form 
$$  \pmatrix
x_1 & x_2 & \!\ldots\! & x_{i-1} & 1 & 0 & \!\ldots\! & 0 & 0
&0 & \!\ldots\! & 0 & 0 & 0 & \!\ldots\! & 0 \cr
y_1 & y_2 & \!\ldots\! & y_{i-1} & 0 & y_{i+1} & \!\ldots\! & y_{j-1} & 1  
&0 & \!\ldots\! & 0 & 0 & 0 & \!\ldots\! & 0 \cr
z_1 & z_2 & \!\ldots\! & z_{i-1} & 0 & z_{i+1} & \!\ldots\! & z_{j-1} & 0  
&z_{j+1} & \!\ldots\! & z_{k-1} & 1 & 0 & \!\ldots\! & 0
\endpmatrix
. \eqno(3)
$$
That is, a matrix with pivots $x_i=y_j=z_k=1$, $i<j<k$, 
$y^i=z^i=z^j=0$, and with no nonzero
entries to the right of the pivots.  Each oriented 3-plane 
corresponds to a unique matrix of this form, or to minus such a matrix.  
For fixed $i,j,k$ we denote
the set of matrices of this type as $e_+(i,j,k)$, and the set of 
negatives of 
these matrices as $e_-(i,j,k)$.  The closures of the sets 
$e_\pm(i,j,k)$, called
Schubert cycles, give a cellular decomposition of $G_N$.

The cell $e_+(i,j,k)$ has dimension $i+j+k-6$.
We give it the orientation \hfill \break
$dx^1 \cdots dx^{i-1} dy^1 \cdots dy^{j-1}
dz^1 \cdots dz^{k-1}$, where of course the variables $y^i, z^i, z^j$ are
skipped in this list.  We orient $e_-(i,j,k)$ so the map $-1 :
e_\pm(i,j,k) \to e_\mp(i,j,k)$ is orientation-preserving.  The
boundary map is then
$$ \eqalign{ 
\partial 
e_\pm(i,j,k) = & (-1)^i e_\pm(i-1,j,k) - e_\mp(i-1,j,k) \cr
&  +  (-1)^{i+j+1} e_\pm(i,j-1,k) + (-1)^{i} e_\mp(i,j-1,k) \cr
&  +  (-1)^{i+j+k+1} e_\pm(i,j,k-1) + (-1)^{i+j} e_\mp(i,j,k-1)} 
\eqno(4) $$
This formula is of course independent of $N$.

$H_4(G_N)$ is then easily computed.  It is $\zz$, and is generated by
$S_N= e_+(1,4,5) + e_+(1,3,6) - e_+(1,2,7)$.  The cycle $\rho(\pi(\nu_N))$ 
doesn't intersect $e_+(1,3,6)$ or $ e_+(1,2,7)$, and hits $e_+(1,4,5)$
at exactly one point, namely
$$ \pmatrix 1 & 0 & 0 & 0 & 0 & 0 & \ldots & 0 \cr
0 & 0 & 0 & 1 & 0 & 0 & \ldots & 0 \cr
0 & 0 & 0 & 0 & 1 & 0 & \ldots & 0
\endpmatrix
, \eqno(5)
$$
and the intersection is transverse.  Thus $\rho(\pi(\nu_N))$ is a
generator of $H^4(G_N)$.  Pulling back we get that $\pi(\nu_N)$ is
a generator of $H^4(B_N)$.  All that remains is to fix the orientation
such that $\pi(\nu_N)$ represents $p_1$.

To fix the orientation we consider the natural embedding $
i: G_N \to G^\comp_{3,N}$, where $G^\comp_{3,N}$ is the Grassmannian
of complex 3-planes in $\comp^N$. The Pontyagin classes on $G_N$ are
pullbacks of Chern classes on $G^\comp_{3,N}$.  
In particular, $p_1 = -i^* c_2$ [MS].  We therefore have only to compute 
the intersection number in $G^\comp_{3,N}$ of $i(S_N)$ with a cycle 
representing $c_2$.  If $W$ is a complex codimension-2 subspace of $\comp^N$, 
then $c_2$ is represented by $Y \subset G^\comp_{3,N}$, the set of 
3-planes in $\comp^N$ whose intersections with $W$ have (complex) dimension
2 or greater [GH].  

If $w_1,\ldots w_N$ are the natural coordinates on $\comp^N$, we choose
$W=\{ w_1 + i w_4=w_2+iw_3 = 0\}$.  A 3-plane spanned by the rows of
$$ \pmatrix x_1 & x_2 & x_3 & x_4 & \ldots \cr
y_1 & y_2 & y_3 & y_4 & \ldots \cr
z_1 & z_2 & z_3 & z_4 & \ldots 
\endpmatrix , \eqno(6)
$$
is in $Y$ if and only if the complex 3-vectors 
$(x_1+ix_4,\; y_1+iy_4, \; z_1+iz_4)$ and
$(x_2+ix_3,\; y_2+iy_3,\; z_2+iz_3)$ are (complex) colinear.  
This is never the case in the closures of $e_+(1,3,6)$ or $ e_+(1,2,7)$. 

Matrices in $e_+(1,4,5)$ take the form
$$ \pmatrix
1 & 0 & 0 & 0 & 0 & 0 & \ldots & 0 \cr
0 & y_2 & y_3 & 1 & 0 & 0 & \ldots & 0 \cr
0 & z_2 & z_3 & 0 & 1 & 0 & \ldots & 0
\endpmatrix , \eqno(7)
$$
$Y$ intersects $e_+(1,4,5)$ at the single point $y_2=y_3=z_2=z_3=0$,
and the intersection number is easily computed to be $+1$.

Thus for a cycle on $G_N$ (or $B_N$) to represent $p_1$, it must be
oriented to intersect $S_N$ (or its image under the natural section)
negatively.  This completes the proof of theorem 2.

\medskip

\noindent {\bf 3. Evaluation of $\mu(p)$.}

\medskip

The finite-dimensional results of section 2
cannot be directly applied to gauge theory.
We need to extend them to appropriate infinite-dimensional
spaces.  Let $H$ be an infinite-dimensional Banach space.  Pick an 
infinite sequence
of linearly independent vectors in $H$.  Then there are natural inclusions
$$ \real^N \inc \real^{N+1} \inc \cdots \inc \real^\infty \inc H,
\eqno(8) $$ 
where $\real^\infty$ is the direct limit of the spaces $\real^N$.
This induces a sequence of inclusions
$$ V_N \inc V_{N+1} \inc \cdots \inc V_\infty \inc V_H \eqno(9)$$
and corresponding inclusions for $V^0$, $B$ and $G$. 
For $N$ large, these inclusions induce isomorphisms in $H_4$ (see
e.g. [MS]), sending $S_N$ to $S_{N+1}$ to $\ldots$ to $S_\infty$ to
$S_H$.
$\pi(\nu_\infty)$ is closed and 
intersects $S_\infty$ once, and $\pi(\nu_H)$ is closed and intersects 
$S_H$ once.  By the same argument as before, we have

\noindent {\bf Theorem 3} {\it  
$\pi(\nu_\infty)$, oriented so as to intersect $S_\infty$
negatively, represents $p_1$ of the bundle $V_\infty \to
B_\infty$, and $\pi(\nu_H)$, oriented to intersect $S_H$ negatively,
represents $p_1$ of $V_H \to B_H$. }

An equivalent description of $p_1$ is as follows.  Let $W$ be a 
codimension-3 subspace in $H$.  Let $Y_W$ be the set of 3-frames 
whose span, intersected with $W$,  is at least 2-dimensional.  When $W=
\{x_1=x_2=x_3=0\}$, $Y_W$ is the same as $\nu_H$.  But, 
since $G_{H^*}$ is connected, 
the choice of $W$ cannot affect the topology of $Y_W$.  Thus $Y_W$, 
oriented to intersect $S_H$ negatively, represents $p_1$ for any 
choice of $W$.

We are now able to construct $\mu$ of a point.  Let $p$ be a point on the
manifold $X$, let $D$ be a geodesic ball around $p$, let $\cala_D$ be the 
$SU(2)$ (or $SO(3)$) connections on $D$ within the Sobolev space $L^q_k$
(the choice of $q$ and $k$ is not important), let $\calg^0$ be the gauge 
transformations in $L^q_{k+1}$ that leave the fiber at $p$ fixed, and let
$\calg$ be all gauge transformations in $L^q_{k+1}$.  Define $\mu_D(p)$ 
to be $-{1 \over 4}p_1$ of the $SO(3)$ bundle $\cala_D/\calg^0
\to \cala_D/\calg$.   $\mu(p)$ is the pullback of $\mu_D(p)$ to $B(X)$ 
via the map that restricts connections on a bundle over $X$ to a bundle
over $D$.

The space $\cala_D/\calg^0$ is isomorphic to the set of connections 
in radial gauge with respect to the point $p$.  In such a gauge the
connection form $A$ vanishes in the radial direction but is otherwise
unconstrained.  In particular, $A(p)=0$, so the curvature at $p$, 
$F_A(p)=dA(p)+ A(p) \wedge A(p) = dA(p)$,
is a linear function of $A$.

Let $H$ be the space of (scalar valued) 1-forms with no radial component.
A connection in radial gauge is 
defined by a triple of elements of $H$, one for each direction in the Lie
Algebra.  Deleting the infinite-codimension
set for which these elements are linearly dependent we get $V_H$. Thus
$\mu_D(p)$ is $-1/4 p_1$ of $V_H \to B_H$, which we have already
computed.  Let $W= \{ \alpha \in H | d^-\alpha(p) = 0 \}$.  Thus $Y_W$ is
the set of connections over $D$, in radial gauge, 
for which the three components
of $F_A^-(p)$ span a 1 (or 0) dimensional subspace of the Lie algebra.
In other words, for which $F_A^-(p)$ is reducible.  Pulling $\mu_D(p)$
back by the restriction map we get the connections on $X$ for which
$F_A^-(p)$ is reducible, i.e. $\nu_p$.  This completes the proof of theorem 1.

\medskip

\noindent{\bf 4. Transversality and Extension to the Boundary.}

\medskip

We have shown that for any point $p$ in our manifold, the cycle 
$\nu_p$ is Poincare dual to $p_1$ of the base point fibration, 
as a class in $\calb^*(X)$.
However, to do Donaldson theory we need more than this.  
Ideally, we want 
$\nu_p$ to intersect the moduli space $\calm_k$ transversely and to 
extend in a well-behaved way to the compactification of moduli space.  
Had we chosen $\nu_p$ 
to depend on $F_A^+$ rather than $F_A^-$, it would still have been 
dual to $p_1$, but would have been useless as a geometric representative 
of $-4\mu(p)$, insofar as $F_A^+$ is identically zero on $\calm_k$.

Even with our definition of $\nu_p$, it is unrealistic to expect
$\nu_p$ to intersect $\calm_k$ transversely for all points $p$.  
For example, if $\calm_k$ has dimension $d<4$, then transversality
would imply  that $\nu_p \cap \calm_k = \emptyset$.  However, there
is a $d+4$ dimensional set of pairs $(A,p)$ for which $F_A^-(p)$
might be reducible.  Since reducibility is a codimension-4 condition,
we should expect reducibility at a $d$-dimensional set of pairs.
Thus for $p$ in a $d$-dimensional subset of $X$, $\nu_p$
would not intersect $\calm_k$ transversely.  There is no reason 
to suppose that this $d$-dimensional set is always empty.

The most we can reasonably expect is the following:

\noindent {\bf Conjecture:} {\it Pick $k>0$ and 
a generic metric on $X$, and let $\calm_k'$ be either  
$\calm_k$ cut down by
standard Donaldson varieties, or $\calm_k$ itself. 
Then, for generic points $p$,
the intersection of $\nu_p$ with $\calm_k'$ is transverse.}

Should this conjecture prove true, then non-transverse intersection
points (for generic metrics) can always be resolved by moving $p$.
If the conjecture is not true, then we will require more subtle 
means of perturbing $\nu_p$, $\calm_k$, or the other Donaldson varieties.   
For many purposes, one wishes
to perturb $\calm_k$ anyway (e.g. modeling connections near the
ends of $\calm_k$ as 
$m$ concentrated charges glued by a particular formula to
connections in $\calm_{k-m}$).  For such purposes, 
the utility of the representative $\nu_p$ does not depend on
the conjecture.

Next we consider the extension of $\nu_p$ to the compactification of 
$\calm_k$.  The boundary of $\calm_k$ consists of strata where $m$ 
instantons have pinched off, leaving a solution of charge $k-m$ behind.  
These take the form $\calm_{k-m} \times S^m(X)$, where $m>0$.  These boundary 
strata have lower dimension than $\calm_k$, so they {\it should} 
not contribute to Donaldson invariants.  To ensure that they do not
contribute, $\nu_p$ must remain a codimension-4 set on
the boundary.  

\noindent {\bf Theorem 4:} {\it The intersection of the closure of
$\nu_p$ with the $m$-th stratum of $\partial \calm_k$ is contained in the 
union of 
$(\nu_p \cap \calm_{k-m}) \times S^m(X)$ and
$\calm_{k-m} \times \{ p \} \times S^{m-1}(X)$.}

\noindent Proof: Consider a sequence of connections $[A_i] \in \calm_k \cap 
\nu_p$  converging to $[A'] \times \{x_1,\ldots,x_m\}$, where $[A'] \in 
\calm_{k-m}$.  If  $p \not \in \{x_i\}$, then $F^-_{A_i}(p)$ converges, 
after suitable gauge transformations, to $F^-_{A'}(p)$.  Since the set 
of rank $\le 1$ matrices is closed and invariant under left 
multiplication by $SO(3)$ (i.e. gauge transformations), $F^-_{A'}$ 
has rank at most 1, and we have the first set.  
If $p \in \{x_i\}$ we are in the second set. QED.   

The first set is manifestly codimension-4.  
If the conjecture holds, 
then, for $m<k$ and generic $p$, the second set
is codimension-4 as well.  
What remains is to consider the first set for $m=k$.  This poses
two difficulties.  First, $\calm_0$
contains the trivial connection (and other reducible connections
if $H_1(X)\ne 0$), and so is not contained in $\calb^*_0$.  
This complication is independent of the choice of representative 
of $\mu(x)$ and is not discussed here.  

(The existence of the trivial connection is also the reason that, 
for $SU(2)$ theory, Donaldson invariants are only well defined for 
$k$ sufficiently large, in the ``stable range''.  For $SO(3)$ theory 
with nontrivial $w_2$, $\calm_0$ is empty, and this complication disappears.) 

The second complication is that every flat connection is in
$\nu_p$, so that $\nu_p$ cannot possibly intersect $\calm_0$
transversely.  To resolve this we must perturb $\calm_0$. If
$\pi_1=0$, so that $\calm_0$ is just the trivial connection, this
is easy.  We just add a small connection that is zero outside a
small neighborhood of $p$.  One can always find a connection for
which $F_A^-(p)$ will be irreducible, so $\nu_p$ will miss the perturbed
$\calm_0$ entirely. If $\pi_1 \ne 0$ and $\calm_0$ contains a 
representation variety of dimension 4 or greater, it may happen that 
one cannot lift $\calm_0$ entirely off $\nu_p$.  In that case
we must interpret ``$\calm_0 \cap \nu_p$'' as the intersection
points that remain after a fixed (but generic) infinitesimal
perturbation of $\calm_0$.

Finally, we consider what must be done if the conjecture
fails.  In that case we would need to construct perturbations $\calm'_k$ of
the moduli spaces $\calm_k$ such that each $\calm'_k$ intersects 
$\nu_p$ transversely, and such that the boundary of $\calm'_k$ 
consists of strata $\calm'_{k-m} \times S^m(X)$.  An analog of 
theorem 4, for $\calm'$, would then follow, and the 
discussion following theorem 4 would also apply.

We close with a sketch of a topological application of this geometric
representative. 
The Donaldson invariants of all known orientable 4-manifolds 
with $b_+>1$ satisfy a recursion relation called
``simple type''. This 
relation roughly says that, given two points $p$ and $q$, $\calm_k\cap\nu_p
\cap\nu_q$ has the same fundamental class as $64 \calm_{k-1}$.
For $p$ and $q$ close and $A$ in $\calm_{k-1}$, one can count the ways to
glue in a concentrated instanton near $p$ and $q$ so as to make the curvature
at $p$ and $q$ reducible.  This number is well short of 64, indicating
that simple type is not just a property of the ends of $\calm_k$, but
involves the topology of the interior as well.
The results of this investigation will appear elsewhere [GS].

I thank Stefan Cordes, Dan Freed, David Groisser, Takashi Kimura, 
Rob Kusner, Tom Parker, Jan Segert, Cliff Taubes and Karen Uhlenbeck for 
extremely helpful discussions.   Part of this work was done at the
1994 Park City/IAS Mathematics Institute.  This work is partially
supported by an NSF Mathematical Sciences Postdoctoral Fellowship
and by Texas Advanced Research Program grant ARP-037.

\smallskip

\refstyle{A}

\enddocument
\bye